\documentclass[12pt,twoside]{extarticle}
\usepackage{mathtext}
\usepackage{titlesec}
\usepackage[cp1251]{inputenc}
\usepackage[T2A]{fontenc}
\usepackage[russian,english]{babel}
\usepackage{fancyhdr}
\usepackage{amsthm}
\usepackage{dsfont}
\usepackage{indentfirst}
\usepackage{mathrsfs}
\usepackage{array}
\usepackage{amssymb}
\usepackage{amsmath}
\usepackage{eucal}
\usepackage{dubloper}
\usepackage{bm}
\usepackage{epsfig}
\usepackage{wrapfig}

\titleformat{\section}[block]{\large\bfseries\filcenter}{\large\bfseries\thesection. }{0pt}{}
\titleformat{\subsection}[block]{\bfseries\filcenter}{\bfseries\thesubsection. }{0pt}{}
\titleformat{\subsubsection}[block]{\bfseries\filcenter}{\bfseries\thesubsubsection. }{0pt}{}
\newcounter{mytheorem}[section]

\edef\lim{\lim\limits}
\renewcommand{\frac}[2]{\dfrac{#1}{#2}}
\numberwithin{equation}{section}
\numberwithin{figure}{section}
\begin{document}

\setlength{\abovedisplayskip}{2mm} % отступ перед выделенной формулой
\setlength{\belowdisplayskip}{2mm} % отступ после выделенной формулой
\parindent=8,8mm%9mm
\renewcommand{\headrulewidth}{0pt}
\edef\sum{\sum\limits}

\thispagestyle{empty}

\begin{center}

\

\

\

\

\textbf{\Large On $\mathrm{K}$-$\mathbb{P}_{t}$-subnormal subgroups of finite groups}

\textbf{\Large and related formations}

\bigskip

\textbf{\bf A.\,F.\,Vasil'ev, T.\,I.\,Vasil'eva}

\end{center}

\bigskip

\begin{center}
{\bf Abstract}
\end{center}

\medskip

\noindent\begin{tabular}{m{8mm}m{138mm}}

&{\small~~~
Let $t$ be a fixed natural number.
A subgroup $H$ of a group $G$ will be called $\mathrm{K}$-$\mathbb{P}_{t}$-subnormal in  $G$ if there exists a chain of subgroups $H = H_{0} \leq H_{1} \leq \cdots \leq H_{m-1} \leq H_{m} = G$ such that either $H_{i-1}$ is normal in $H_{i}$ or $|H_{i} : H_{i-1}|$ is a some prime $p$ and $p-1$ is not divisible by the $(t+1)$th powers of primes for every $i = 1,\ldots , n$.
In this work, properties of $\mathrm{K}$-$\mathbb{P}_{t}$-subnormal subgroups and classes of groups with Sylow $\mathrm{K}$-$\mathbb{P}_{t}$-subnormal subgroups are obtained.

\medskip

}

\medskip

\smallskip

MSC2010\  20D10, 20E15, 20F16

%20D10: Solvable groups, theory of formations, Schunck classes, Fitting classes, $\pi$-length, ranks [See also 20F17]
%20D20 Sylow subgroups, Sylow properties, -groups, -structure
%20D25 Special subgroups (Frattini, Fitting, etc.)
%20D30 Series and lattices of subgroups
%20D35 Subnormal subgroups
%20D40 Products of subgroups
%20E15: Chains and lattices of subgroups, subnormal subgroups
%20F16: Solvable groups, supersolvable groups [See also 20D10]
%20F17: Formations of groups, Fitting classes [See also 20D10]

\end{tabular}

\bigskip

{\parindent=0mm
\textbf{\bf {\large Introduction}}
}
\bigskip

All groups under consideration are finite.

The goal of this work is to study the following generalization of the concept of a subnormal subgroup in a group and to find some of its applications.

\medskip

{\textbf{Definition 1.}
{\it Let $t$ be a fixed natural number.
A subgroup $H$ of a group $G$ will be called $\mathrm{K}$-$\mathbb{P}_{t}$-subnormal in  $G$ if there exists a chain of subgroups

\smallskip

\ \ \ \ \ \ \ \ \ \ \ \ \ \ \ \ $H = H_{0} \leq H_{1} \leq \cdots \leq H_{n-1} \leq H_{n} = G$ \hspace{\stretch{1}}$(1.1)$

\smallskip

{\parindent=0mm
such that either $H_{i-1}\unlhd H_{i}$ or $|H_{i} : H_{i-1}|$ is a some prime $p$ and $p-1$ is not divisible by the $(t+1)$th powers of primes for every $i = 1,\ldots , n$.
%We say that $H$ is a $\mathrm{K}$-$\mathbb{P}_{t}$-$sn$ $G$.
}
}
\medskip

A subgroup $H$ of a group $G$ is said to be:

{\it $\mathbb{P}$-subnormal in $G$} \cite{VasVasTyu2010} if either $H= G$ or there exists a chain of subgroups $(1.1)$
such that $|H_{i} : H_{i-1}|$ is a prime for every $i = 1,\ldots , n$;

{\it $\mathrm{K}$-$\mathbb{P}$-subnormal in $G$} \cite{VasVasTyu2014} if there exists a chain of subgroups $(1.1)$
such that either $H_{i-1}\unlhd H_{i}$ or $|H_{i} : H_{i-1}|$ is a prime for every $i = 1,\ldots , n$.

\medskip

In recent years, the concepts $\mathbb{P}$-subnormal and $\mathrm{K}$-$\mathbb{P}$-subnormal
subgroups have been used in the works of many authors to solve various problems in group theory (see, for example, [3--17]).
%\cite{KniMon13, VasVasMys16, VasVasPar18, BalLiPedSU20, Mur20, LucNem21, VTI22, BalMadShuPed22, BalMadPedWU23, CheZhaLi23, VTIKor23, Mur24, Lis24}).

It is clear that every $\mathrm{K}$-$\mathbb{P}_{t}$-subnormal subgroup in $G$ is $\mathrm{K}$-$\mathbb{P}$-subnormal. The converse statement does not hold in the general case. For example, let $t=3$ and let $G=AB$ be a group, where $A\cong Z_{17}$ and $B\cong Aut(Z_{17})\cong Z_{16}$. In $G$, there is a chain of subgroups
$$H \unlhd B < G,$$
with $|H|=2$. Then $H$ is $\mathrm{K}$-$\mathbb{P}$-subnormal and $\mathbb{P}$-subnormal in $G$, but $H$ is not $\mathrm{K}$-$\mathbb{P}_{3}$-subnormal in $G$, since $|G:B|=17$, $17-1=2^{4}$.

\medskip

Let $\mathfrak{F}$ be a non-empty formation. A subgroup $H$ of a group $G$ is said to be:

{\it $\mathfrak{F}$-subnormal in $G$} \cite{BalEzq} if either $H=G$ or there exists a chain of subgroups $(1.1)$
such that $H_{i}^{\mathfrak{F}} \leq H_{i-1}$ for every $i = 1,\ldots , n$.
Here $G^{\mathfrak{F}}$ is the $\mathfrak{F}$-residual of $G$, i.e. the least normal subgroup of $G$ for which $G/G^{\mathfrak{F}}\in \mathfrak{F}$;

{\it $\mathrm{K}$-$\mathfrak{F}$-subnormal in $G$} \cite{BalEzq} if there exists a chain of subgroups $(1.1)$
such that either $H_{i-1}\unlhd H_{i}$ or $H_{i}^{\mathfrak{F}} \leq H_{i-1}$ for every $i = 1,\ldots , n$.

\medskip

Let $\mathfrak{U}_{k}$ be the class of all supersoluble groups in which exponents are not divided by the $(k + 1)$th powers of primes, where $k$ is a naturel number. In \cite{MonSok_Sib_el23}, V.S. Monakhov and I.L. Sochor showed that $\mathfrak{U}_{k}$ is a hereditary formation, and studied the class of groups $\mathrm{w}\mathfrak{U}_{k}$ in which every Sylow subgroup is $\mathfrak{U}_{k}$-subnormal.

In a soluble group $G$ every $\mathrm{K}$-$\mathbb{P}_{t}$-subnormal subgroup is $\mathfrak{U}_{t} $-subnormal in $ G$ (Lemma 2.4 ). The converse does not hold (Example 2.1).

Note that in $G$, a $\mathrm{K}$-$\mathbb{P}_{t}$-subnormal subgroup is not $\mathfrak{U}_{t}$-subnormal in general case.
For example, let $t=2$ and let $G\cong A_{5}$ is an alternating group of degree 5. The Sylow 2-subgroup $H$ of $G$ is $\mathrm{K}$-$\mathbb{P}_{2}$-subnormal, but is not $\mathfrak{U}_{2}$-subnormal in $G$, since $H\unlhd H_{1} < G$, where $H_{1}\cong A_{4} $, $5-1=2^{2}$, $H_{1}^{\mathfrak{U}_{2}} = H$, $G^{\mathfrak{U}_ {2}} = G$.

In this work, properties of $\mathrm{K}$-$\mathbb{P}_{t}$-subnormal subgroups and classes of groups with Sylow $\mathrm{K}$-$\mathbb{P}_{t}$-subnormal subgroups are obtained.

%Note that every $\mathrm{K}$-$\mathbb{P}_{t}$-subnormal subgroup in $G$ is $\mathrm{K}$-$\mathbb{P}$-subnormal. By \cite{VasVasTyu2014}, a subgroup $H$ of a group $G$ is called $\mathrm{K}$-$\mathbb{P}$-subnormal in  $G$ if there exists a chain of subgroups $(1.1)$ such that either $H_{i-1}\unlhd H_{i}$ or $|H_{i} : H_{i-1}|$ is a prime number for every $i = 1,\ldots , n$.

\bigskip

{\parindent=0mm
\textbf{\bf {\large 1. Preliminary results}}
}

\bigskip

We use the notation and terminology from \cite{BalEzq, DH}.
We recall some concepts significant in the paper.

Let $G$ be a group. If $H$ is a subgroup of $G$, we write $H\leq G$ and if $H\not= G$,
we write $H < G$.
We denote by $|G|$ the order of $G$;
by $\pi(G)$ the set of all distinct prime divisors of the order of $G$;
by ${\rm {Syl}}_p(G)$ the set of all Sylow $p$-subgroups of $G$;
by ${\rm {Syl}}(G)$ the set of all Sylow subgroups of $G$;
by $\mathrm{Core}_G(M)$ the core of subgroup $M$ in $G$, i.e. $\mathrm{Core}_G(M)=\cap M^{x}$ for all $x\in G$;
%intersection of all subgroups conjugated with $M$ in $G$;
by $F(G)$ the Fitting subgroup of $G$;
by $|G:H|$ the index of $H$ in $G$;
by $\pi(G:H)$ the set of all different prime divisors of $|G:H|$;
%For a subgroup $M$ of a group $G$, $M_{G}$ denotes the core of $M$ in $G$, i.e. $M_{G}=\cap M^{x}$ for all $x\in G$;
%произведение всех нормальных ниль\-по\-тент\-ных подгрупп группы $G$;
%by $F_p(G)$ the $p$-nilpotent radical of $G$, i.e.
%the product of all normal $p$-nilpotent subgroups of $G$.
by $Z_n$ the cyclic group of order $n$;
by $\mathbb{P}$ the set of all primes;
%$\pi$ denotes a set of primes;
%$\pi'=\Bbb{P}\setminus\pi$;
by $\mathfrak{S}$ the class of all soluble groups;
by $\mathfrak{U}$ the class of all supersoluble groups;
by $\mathfrak{N}$ the class of all nilpotent groups;
by $\mathfrak{N}_{p}$ the class of all $p$-groups for $p\in\mathbb{P}$;
by $\mathfrak{A}(p-1)$ the class of all abelian groups of exponent dividing $p-1$.

A group $G$ of order $p_1^{\alpha_1}p_2^{\alpha_2}\cdots p_n^{\alpha_n}$
(where $p_i$ is a prime, $i=1, 2,\ldots, n$) is called {\it Ore dispersive} %\cite[p.~251]{Shem},
or {\it Sylow tower group} whenever $p_1>p_2>\cdots >p_n$ and
$G$ has a normal subgroup of order $p_1^{\alpha_1}p_2^{\alpha_2}\cdots p_i^{\alpha_i}$
for every $i=1, 2,\ldots, n$.

A class of groups $\mathfrak{F}$ is called a {\it formation} if the following conditions hold:
(a) every quotient group of a group lying in $\mathfrak{F}$ also lies in
$\mathfrak{F}$; (b) if $G/N_{i}\in \mathfrak{F}$, $N_{i}\unlhd G$, $i=1, 2$ then $G/N_{1}\cap N_{2}\in \mathfrak{F}$.
A formation $\mathfrak{F}$ is called
{\it hereditary} whenever $\mathfrak{F}$ together with every group contains
all its subgroups, and {\it saturated}, if $G/\Phi(G)\in\mathfrak{F}$ implies that $G\in\mathfrak{F}$.

A function $f:\Bbb{P}\rightarrow \{$formations$\}$ is called a {\it local function}.
A formation $\mathfrak{F}$ is called {\it local},
if there exists a local function $f$ such that $\mathfrak{F}=LF(f)=(G | G/C_{G}(H/K)\in f(p)$ for every chief factor $H/K$ of a group $G$ and for all primes $p\in\pi(H/K))$.

\medskip

%{\bf Lemma~1.1 \cite[Lemma~3.9]{Shem}.} {\it If $H/K$ is a chief factor of a group $G$ and $p\in\pi(H/K)$ then $G/C_G(H/K)$ doesn't contain nonidentity normal $p$-subgroup and besides $F_p(G)\leq C_G(H/K)$.}

\medskip

{\bf Lemma 1.1} \cite[Chap. I, Theorem~1.4]{Wei}.
{\it Let  $H/K$ be a chief $p$-factor of a group $G$. $|H/K|=p$ if and only if
${\rm Aut}_G(H/K)$ is abelian group of exponent dividing $p-1$.
}

\medskip

%We will need the following known results.

{\bf Lemma 1.2} \cite[Lemma 3.6, Corollary 4.3.1]{VasVasTyu2014}.
{\it The class $\overline{\mathrm{w}}\mathfrak{U}$ of all groups in which any Sylow subgroup is $\mathrm{K}$-$\mathbb{P}$-subnormal consists of Ore dispersive groups, and forms a hereditary saturated formation.
}

\medskip

{\bf Lemma 1.3 } \cite[Theorem 2]{Vas2004}.
{\it A group $G$ is supersoluble if and only if $G$ can be represented as the product of two nilpotent $\mathbb{P}$-subnormal subgroups.
}

\medskip

We will need some properties of $\mathfrak{F}$-subnormal subgroups (see, for example, \cite[Chap. 6]{BalEzq}). In what follows, $\mathfrak{F}$ means a non-empty formation.

{\bf Lemma 1.4.}
{\it Let $H$ and $U$ be subgroups of a group $G$ and let $N\unlhd G$. Suppose that $\mathfrak{F}$ is a formation. Then the following statements hold.

$(1)$ If $H$ is $\mathfrak{F}$-subnormal in $U$ and $U$ is $\mathfrak{F}$-subnormal in $G$ then $H$ is $\mathfrak{F}$-subnormal in $G$.

$(2)$ If $N \leq U$ and $U/N$ is $\mathfrak{F}$-subnormal in $G/N$ then $U$ is $\mathfrak{F}$-subnormal in $G$.

$(3)$ If $H$ is $\mathfrak{F}$-subnormal in $G$ then $HN/N$ is $\mathfrak{F}$-subnormal in $G/N$.

Suppose that $\mathfrak{F}$ is a a hereditary formation. Then the following statements hold.

$(4)$ If $G^{\mathfrak{F}} \leq H$ then $H$ is $\mathfrak{F}$-subnormal in $G$.

$(5)$ If $H$ is $\mathfrak{F}$-subnormal in $G$ then $H\cap U$ is $\mathfrak{F}$-subnormal in $U$.
}

\medskip

According to \cite{VasVas11}, class of groups $\mathrm{w}\mathfrak{F} = (G | \pi(G)\subseteq \pi(\mathfrak{F})$ and every Sylow subgroup of $G$ is  $\mathfrak{F}$-subnormal in $G)$. Here $\pi(\mathfrak{F})$ is the set of all distinct prime divisors $|G|$ for $G\in\mathfrak{F}$.

\medskip

{\bf Lemma 1.5} \cite[Lemma 1.6]{VasVas11}.
{\it If $\mathfrak{F}$ is a hereditary formation and $\mathfrak{F}\subseteq \mathfrak{S}$ then  $\mathrm{w}\mathfrak{F}\subseteq \mathfrak{S}$.
}

\medskip

{\bf Lemma 1.6} \cite[Theorem B]{VasVas11}.
{\it If $\mathfrak{F}$ is a hereditary saturated formation then $\mathrm{w}\mathfrak{F}$ is a hereditary saturated formation.
}

\bigskip

{\parindent=0mm
\textbf{\bf {\large 2. Properties of $\mathrm{K}$-$\mathbb{P}_{t}$-subnormal subgroups of groups}}
}

\bigskip

Further $t$ means a fixed natural number.

\medskip

{\bf Lemma 2.1.}
{\it Let $H$ be a subgroup of a group $G$. Then the following statements hold.

$(1)$ If $H$ is $\mathrm{K}$-$\mathbb{P}_{t}$-subnormal in $G$, then $H^{x}$ is $\mathrm{K}$-$\mathbb{P}_{t}$-subnormal in $G$ for all $x\in G$.

$(2)$ If $H\leq R\leq G$, $H$ is $\mathrm{K}$-$\mathbb{P}_{t}$-subnormal in $R$ and $R$ is $\mathrm{K}$-$\mathbb{P}_{t}$-subnormal in $G$, then $H$ is $\mathrm{K}$-$\mathbb{P}_{t}$-subnormal in $G$.
}

\medskip

{\it Proof.}
%\begin{proof}
Statements (1) and (2) follow from Definition 1.
\hspace{\stretch{1}}$\square$
%{\qedsymbol}
%\end{proof}

\medskip

{\bf Lemma 2.2.}
{\it Let $H$ be a subgroup of a group $G$ and $N\unlhd G$. Then the following statements hold.

$(1)$ If $H$ is $\mathrm{K}$-$\mathbb{P}_{t}$-subnormal in $G$, then $(H\cap N)$ is $\mathrm{K}$-$\mathbb{P}_{t}$-subnormal in $N$ and $HN/N$ is $\mathrm{K}$-$\mathbb{P}_{t}$-subnormal in $G/N$.

$(2)$ If $N\leq H$ and $H/N$ is $\mathrm{K}$-$\mathbb{P}_{t}$-subnormal in $G/N$, then $H$ is $\mathrm{K}$-$\mathbb{P}_{t}$-subnormal in $G$.

$(3)$ The subgroup $HN$ is $\mathrm{K}$-$\mathbb{P}_{t}$-subnormal in $G$ if and only if $HN/N$ is $\mathrm{K}$-$\mathbb{P}_{t}$-subnormal in $G/N$.}

$(4)$  If $HN_{i}$ is $\mathrm{K}$-$\mathbb{P}_{t}$-subnormal in $G$ and $N_{i}\unlhd G$, $i = 1, 2$, then $(HN_{1}\cap HN_{2})$ is $\mathrm{K}$-$\mathbb{P}_{t}$-subnormal in $G$.
}

\medskip

{\it Proof.}
%\begin{proof}
(1) Suppose that $H$ is $\mathrm{K}$-$\mathbb{P}_{t}$-subnormal in $G$. We can assume that $H\not= G$. There is a chain of subgroups (1.1).
Consider the chains of subgroups
$$H\cap N = H_{0}\cap N \leq H_{1}\cap N \leq \cdots \leq H_{n-1}\cap N \leq H_{n}\cap N = N,$$
$$HN/N = H_{0}N/N \leq H_{1}N/N \leq \cdots \leq H_{n-1}N/N \leq H_{n}N/N = G/N.$$

If $H_{i-1}\cap N \unlhd H_{i}\cap N$ for all $i = 1, \ldots , n$, then $H\cap N$ is $\mathrm{K}$-$\mathbb{P}_{t}$-subnormal in $N$.

Let's assume that $H_{i-1}\cap N\ntrianglelefteq H_{i}\cap N$ for some $i\in\{1, \ldots , n\}$. In that case $H_{i-1}\cap N\not= H_{i}\cap N$ and $H_{i-1}\ntrianglelefteq H_{i}$. Then $|H_{i} : H_{i-1}|$  is a some prime $p$ and $p-1$ is not divisible by the $(t+1)$th powers of primes. Then $H_{i-1}$ is maximal in $H_{i}$ and $H_{i}=(H_{i}\cap N)H_{i-1}$. We have $|H_{i}\cap N : H_{i-1}\cap N|=|H_{i} : H_{i-1}|=p$. Thus, $H\cap N$ is $\mathrm{K}$-$\mathbb{P}_{t}$-subnormal in $N$.

If $H_{i-1}N/N \unlhd H_{i}N/N$ for all $i = 1, \ldots , n$, then $HN/N$ is $\mathrm{K}$-$\mathbb{P}_{t}$-subnormal in $G/N$.
Let's assume that $H_{j-1}N/N\ntrianglelefteq H_{j}N/N$ for some $j\in\{1, \ldots , n\}$. Then $H_{j-1}\ntrianglelefteq H_{j}$. By Definition 1 $|H_{j} : H_{j-1}|$  is a some prime $q$ and $q-1$ is not divisible by the $(t+1)$th powers of primes. Hence $H_{j-1}$ is maximal in $H_{j}$ and $H_{j-1}=(H_{j}\cap N)H_{j-1}$. We have $|H_{j}N/N : H_{j-1}N/N|=|H_{j} : H_{j-1}|=q$. Therefore $HN/N$ is $\mathrm{K}$-$\mathbb{P}_{t}$-subnormal in $G/N$.

(2) Assume that $H/N$ is $\mathrm{K}$-$\mathbb{P}_{t}$-subnormal in $G/N$. Then there exists a chain of subgroups
$$H/N = H_{0}/N \leq H_{1}/N \leq \cdots \leq H_{n-1}/N \leq H_{n}/N = G/N$$
such that either $H_{i-1}/N\unlhd H_{i}/N$ or $|H_{i}/N : H_{i-1}/N|$ is a some prime $p$ and $p-1$ is not divisible by the $(t+1)$th powers of primes for every $i = 1,\ldots , n$. Therefore either $H_{i-1}\unlhd H_{i}$ or $|H_{i}/N : H_{i-1}/N|=|H_{i} : H_{i-1}|$ and $H$ is $\mathrm{K}$-$\mathbb{P}_{t}$-subnormal in $G$.

(3) The statement follows from statements (1) and (2).

(4) Let $HN_{i}$ is $\mathrm{K}$-$\mathbb{P}_{t}$-subnormal in $G$ and $N_{i}\unlhd G$, $i = 1, 2$. Then there is a chain of subgroups
$$HN_{1} = K_{0} \leq K_{1} \leq \cdots \leq  K_{s-1} \leq  K_{s} = G$$
such that either $K_{i-1}$ is normal in $K_{i}$ or $|K_{i} : K_{i-1}|$ is a some prime $p$ and $p-1$ is not divisible by the $(t+1)$th powers of primes for every $i = 1,... , s$. If
$(K_{i-1}\cap HN_{2}) \unlhd (K_{i} \cap HN_{2})$ for every $i = 1,... , s$,
then $HN_{1} \cap HN_{2}$ is $\mathrm{K}$-$\mathbb{P}_{t}$-subnormal in $G$. Let there be $j \in \{1,... ,s\}$ such that $(K_{j-1}\cap HN_{2}) \ntrianglelefteq (K_{j} \cap HN_{2})$.
Then $K_{j-1} \ntrianglelefteq K_{j}$ and $|K_{j} : K_{j-1}|$ is a prime $p$ and $p-1$ is not divisible by the $(t+1)$th powers of primes. Since
$K_{j-1} \leq (K_{j} \cap HN_{2})K_{j-1} \leq K_{j}$
and $(K_{j} \cap HN_{2})K_{j-1} = (K_{j} \cap N_{2})K_{j-1}$ is a subgroup of $K_{j}$, it follows that $(K_{j} \cap HN_{2})K_{j-1} = K_{j}$.
Therefore,
$|K_{j} \cap HN_{2} : K_{j-1} \cap HN_{2}| = |(K_{j} : K_{j-1}|=p$. Thus, considering the chain of subgroups
$$HN_{1}\cap HN_{2} = K_{0} \cap HN_{2} \leq K_{1} \cap HN_{2} \leq \cdots \leq K_{s-1} \cap HN_{2} \leq K_{s} \cap HN_{2} = HN_{2}$$
and taking into account the $\mathrm{K}$-$\mathbb{P}_{t}$-subnormality of $HN_{2}$ in $G$, we see that $(HN_{1} \cap HN_{2})$ is $\mathrm{K}$-$\mathbb{P}_{t}$-subnormal in $G$.
\hspace{\stretch{1}}$\square$
%{\qedsymbol}
%\end{proof}

\medskip

{\bf Lemma 2.3.}
{\it Let $G$ a soluble group. Let $H$ and $U$ be subgroups of $G$. Then the following statements hold.

$(1)$ If $H$ is $\mathrm{K}$-$\mathbb{P}_{t}$-subnormal in $G$, then $(H\cap U)$ is $\mathrm{K}$-$\mathbb{P}_{t}$-subnormal in $U$.

$(2)$ If $H$ is $\mathrm{K}$-$\mathbb{P}_{t}$-subnormal in $G$ and $U$ is $\mathrm{K}$-$\mathbb{P}_{t}$-subnormal in $G$, then $(H\cap U)$ is $\mathrm{K}$-$\mathbb{P}_{t}$-subnormal in $G$.
}

\medskip

{\it Proof.}
%\begin{proof}
(1) Let $H$ be $\mathrm{K}$-$\mathbb{P}_{t}$-subnormal in $G$.

If $H_{i-1}\cap U \unlhd H_{i}\cap U$ for all $i=1,\ldots, n$ and $H_{i}$ is from the chain (1.1), then $(H\cap U)$ is $\mathrm{K}$-$\mathbb{P}_{t}$-subnormal in $U$.

Suppose that there exists $j\in\{1,\ldots, n\}$ such that $H_{j-1}\cap U\ntrianglelefteq H_{j}\cap U$. Then $H_{j-1}\ntrianglelefteq H_{j}$ and $|H_{j}:H_{j-1}|$ is a some prime $p$ and $p-1$ is not divisible by the $(t+1)$th powers of primes.

Let's denote $L=\mathrm{Core}_{H_{j}}(H_{j-1})$. By \cite[Chap. A, Theorem 15.2(4)]{DH} ${H_{j}}/L=N/L\cdot {H_{j-1}}/L$, where $N/L$ is a minimal normal subgroup of ${H_{j}}/L$, $N/L=C_{{H_{j}}/L}(N/L)$ and $(N/L)\cap ({H_{j-1}}/L)=L/L$. Then $|N/L|=|H_{j}:H_{j-1}|=p$ and ${H_{j-1}}/L\cong {H_{j}}/L/C_{{H_{j}}/L}(N/L)\cong Aut_{H_{j}}(N/L)$ is isomorphic to a subgroup from $Z_{p-1}$. Therefore $N/L\in \mathrm{Syl}_{p}(G)$.

Let us first assume that $p$ does not divide $|(H_{j}\cap U)L/L|$. Since $G$ is soluble, by Hall's theorem \cite[Chap. A, Theorem 3.3]{DH} $(H_{j}\cap U)L/L\leq (H_{j-1}/L)^{xL}$ for some $x\in H_{j}$. Then $(H_{j}\cap U)/(U\cap L)\cong (H_{j}\cap U)L/L$ is cyclic, we have $(H_{j-1}\cap U)/(U\cap L)\unlhd (H_{j}\cap U)/(U\cap L)$. We have obtained a contradiction with the assumption $H_{j-1}\cap U\ntrianglelefteq H_{j}\cap U$.

Thus $p$ divides $|(H_{j}\cap U)L/L|$. Since $N/L\unlhd {H_{j}}/L$, by Sylow's theorem $N/L\leq (H_{j}\cap U)L/L$. Then
$$(H_{j}\cap U)L/L=(H_{j}\cap U)L/L\cap N/L\cdot {H_{j-1}}/L=N/L\cdot ({H_{j-1}\cap U})L/L$$
and $p=|N/L|=|(H_{j}\cap U):(H_{j-1}\cap U)|$.
Therefore we have the chain of subgroups
$$H\cap U = H_{0}\cap U \leq H_{1}\cap U \leq \cdots \leq H_{n-1}\cap U \leq H_{n}\cap U = U$$
such that either $(H_{i-1}\cap U)\unlhd (H_{i}\cap U)$ or $|(H_{i}\cap U) : (H_{i-1}\cap U)|$ is a some prime $p$ and $p-1$ is not divisible by the $(t+1)$th powers of primes for every $i = 1,\ldots , n$.

Statement (2) follows from (1) and Lemma 2.1(2).
\hspace{\stretch{1}}$\square$
%{\qedsymbol}
%\end{proof}

\medskip

{\bf Lemma 2.4.} {\it If $G$ is a soluble group and $H$ is $\mathrm{K}$-$\mathbb{P}_{t}$-subnormal in $G$, then $H$ is $\mathfrak{U}_{t} $-subnormal in $G$.}

\medskip

{\it Proof.}
%\begin{proof}
By Definition 1 there is a chain of subgroups (1.1). Assume that $H\not= G$. Then $H_{i-1}\not= H_{i}$ for some $i\in\{ 1, \ldots , n\}$.

Assume that  $H_{i-1}\unlhd H_{i}$. Since $G$ is soluble, we draw a composition series through $H_{i-1}$ and $H_{i}$:
$$H_{i-1} = R_{0}< R_{1}< \cdots < R_{m-1} < R_{m} = H_{i},$$
there $R_{j-1}$ is normal in $R_{j}$ and $|R_{j}:R_{j-1}|$ is some prime $q$, $j=1,\ldots, m$. Since $R_{j-1}=\mathrm{Core}_{R_{j}}(R_{j-1})$, we have $R_{j}/R_{j-1}\in\mathfrak{U}_{t}$. Thus ${R_{j}}^{\mathfrak{U}_{t}}\leq R_{j-1}$.

Assume that $H_{i-1}\ntrianglelefteq H_{i}$. Then $|H_{i} : H_{i-1}|$ is a some prime $p$ and $p-1$ is not divisible by the $(t+1)$th powers of primes. Since  $L=\mathrm{Core}_{H_{i}}(H_{i-1})\not=H_{i-1}$, we have ${H_{i}}/L=N/L\cdot {H_{i-1}}/L$, where $N/L$ is a minimal normal subgroup of ${H_{i}}/L$, $N/L=C_{{H_{i}}/L}(N/L)$ and $(N/L)\cap ({H_{i-1}}/L)=L/L$ by \cite[Chap. A, Theorem 15.2(4)]{DH}. Then $|N/L|=|H_{i}:H_{i-1}|=p$ and ${H_{i-1}}/L\cong {H_{i}}/L/C_{{H_{i}}/L}(N/L)\cong Aut_{H_{i}}(N/L)$ is isomorphic to a subgroup from $Z_{p-1}$. Thus ${H_{i}}/L\in\mathfrak{U}_{t}$ and ${H_{i}}^{\mathfrak{U}_{t}}\leq L\leq H_{i-1}$.

Therefore $H$ is $\mathfrak{U}_{t}$-subnormal in $G$.
\hspace{\stretch{1}}$\square$
%{\qedsymbol}
%\end{proof}

\medskip

Note that the converse of Lemma 2.4 does not always hold.

{\bf Example 2.1.} Let $t=1$ and let $G$  be a non-abelian group of order 39. In $G$, the Sylow 3-subgroup $H$ is $\mathfrak{U}_{1}$-subnormal, but not $\mathrm{K}$-$\mathbb {P}_{1}$ is subnormal, since $G^{\mathfrak{U}_{1}}=1\leq H$, $|G:H|=13$ and $13-1=2^ {2}\cdot 3$.

\bigskip

\textbf{\bf {\large 3. Classes of groups with $\mathrm{K}$-$\mathbb{P}_{t}$-subnormal Sylow subgroups}}

\bigskip

{\bf Lemma 3.1.} {\it  Let $G$ be a group. Let $H\in\mathrm{Syl}_{p}(G)$ and $N\unlhd G$. If $H$ is $\mathrm{K}$-$\mathbb{P}_{t}$-subnormal in $G$, then every Sylow $p$-subgroup of $N$ is $\mathrm{K}$-$\mathbb{P}_{t}$-subnormal in $N$ and every Sylow $p$-subgroup of $G/N$ is $\mathrm{K}$-$\mathbb{P}_{t}$-subnormal in $G/N$.
}

\medskip

{\it Proof.}
%\begin{proof}
Let $P\in\mathrm{Syl}_{p}(N)$ and $R/N\in\mathrm{Syl}_{p}(G/N)$. By Sylow's theorem $P\leq H^{x}$ for some $x\in G$ and $R/N=(HN/N)^{yN}$ for some $y\in G$. Then $P=H^{x}\cap N$ is $\mathrm{K}$-$\mathbb{P}_{t}$-subnormal in $N$ and $R/N=H^{y}N/N$ is $\mathrm{K}$-$\mathbb{P}_{t}$-subnormal in $G/N$ by Lemmas 2.1(1), 2.2(1).
\hspace{\stretch{1}}$\square$
%{\qedsymbol}
%\end{proof}

\medskip

By Lemma 1.2 the class of groups $\mathrm{w}_{t}\mathfrak{U}$ consists of Ore dispersive groups, therefore $\mathrm{w}_{t}\mathfrak{U}\subseteq\mathfrak{S}$.

\medskip

{\bf Theorem 3.1.} {\it  Let $\mathfrak{H}$ be the class of all groups in which every Sylow subgroup is $\mathrm{K}$-$\mathbb{P}_{t}$-subnormal
Let $G$ be a group. Then the following assertions hold.

$(1)$ $\mathfrak{N}\subseteq\mathfrak{H}$ and $\mathfrak{H}$ consists of Ore dispersive groups.

$(2)$ If $G\in\mathfrak{H}$ and $N\unlhd G$, then $G/N\in\mathfrak{H}$.

$(3)$ If $G/N_{1}\in\mathfrak{H}$ and $G/N_{2}\in\mathfrak{H}$ for any $N_{i}\unlhd G$, $i = 1,2$, then $G/N_{1}\cap N_{2}\in\mathfrak{H}$.

$(4)$ A direct product of groups from $\mathfrak{H}$ lies in $\mathfrak{H}$.

$(5)$ If $G\in\mathfrak{H}$ and $U$ is a subgroup of $G$, then $U\in\mathfrak{H}$.

$(6)$ If If $G/\Phi(G)\in\mathfrak{H}$, then $G\in\mathfrak{H}$.

$(7)$ The class $\mathfrak{H}$ is a hereditary saturated formation.
}

\medskip

{\it Proof.}
%\begin{proof}
Statement (1) follows from Definition 1 and Lemma 1.2.

Statement (2) follows from Lemma 3.1.

(3) Let $G$ be a group of the least order such that $G/N_{i}\in\mathfrak{H}$, $N_{i}\unlhd G$, $i = 1,2$, but $G/N_{1}\cap N_{2}\not\in\mathfrak{H}$.

If $N = N_{1}\cap N_{2}\not= 1$, then $G/N/N_{i}/N\cong G/N_{i}\in \mathfrak{H}$ for i = 1,2. It follows from
$|G/N| < |G|$ that
$G/N/(N_{1}/N\cap N_{2}/N)\cong G/N_{1}\cap N_{2}\in\mathfrak{H}$. This contradicts the choice of $G$.

Thus, $N_{1}\cap N_{2} = 1$. Let $Q\in\mathrm{Syl}_{q}(G)$. Then $QN_{i}/N_{i}\in\mathrm{Syl}_{q}(G/N_{i})$ and $QN_{i}/N_{i}$ is $\mathrm{K}$-$\mathbb{P}_{t}$-subnormal in $G/N_{i}$ for $i=1, 2$. By Lemma 2.2(2) $QN_{i}$ is $\mathrm{K}$-$\mathbb{P}_{t}$-subnormal in $G$ for $i=1, 2$. Since $G$ is soluble, by Lemma 2.3(2) $QN_{1}\cap QN_{2}$ is $\mathrm{K}$-$\mathbb{P}_{t}$-subnormal in $G$. By \cite[Chap. A, Theorem 6.4(b)]{DH} $QN_{1}\cap QN_{2}=Q(N_{1}\cap N_{2})=Q$ is $\mathrm{K}$-$\mathbb{P}_{t}$-subnormal in $G$. Therefore $G=G/N_{1}\cap N_{2}\in\mathfrak{H}$. The contradiction thus obtained completes the proof of Statement (3).

Statement (4) follows from (3).

(5) Let $G\in\mathfrak{H}$ and $U\leq G$. Let's take $P\in\mathrm{Syl}_{p}(U)$. By Sylow's theorem $P\leq P_{1}$ for some $P_{1}\in\mathrm{Syl}_{p}(G)$. Since $G$ is soluble and $P_{1}$ is $\mathrm{K}$-$\mathbb{P}_{t}$-subnormal in $G$, by Lemma 2.3(1) $P=P_{1}\cap U$ is $\mathrm{K}$-$\mathbb{P}_{t}$-subnormal in $U$. Therefore $U\in\mathfrak{H}$. Statement (5) has been proven.

(6) Let $G$ be a group of the least order for which $G/\Phi(G)\in\mathfrak{H}$ and $G\not\in\mathfrak{H}$. Then $G$ is soluble, since $G/\Phi(G)$ and $\Phi(G)$ are soluble. Let $N$ be a minimal normal subgroup of $G$. We note that $N$ is a $p$-group for some prime $p$. By \cite[Chap. A, Theorem 9.2(e)]{DH}  $\Phi(G)N/N\leq \Phi(G/N)$. Since $G/\Phi(G)N\in\mathfrak{H}$, we have $(G/N)/\Phi(G/N)\in\mathfrak{H}$. From $|G/N| < |G|$, it follows
$G/N\in\mathfrak{H}$.

By (2) and (3) $\mathfrak{H}$ is a formation. Then $N$ is a unique minimal normal subgroup of $G$. Thus, $N\leq \Phi(G)\leq F(G)$ and $F(G)$ is $p$-groups. By \cite[Chap. A, Theorem 10.6(c)]{DH} $\Phi(G) < F(G)$. Let $Q\in\mathrm{Syl}_{q}(G)$. Since $QN/N\in\mathrm{Syl}_{q}(G/N)$ and $QN/N\nleqslant\Phi(G/N)$, we have $QN/N\not=N/N$ and $QN/N$ is $\mathrm{K}$-$\mathbb{P}_{t}$-subnormal in $G/N$. By Lemma 2.2(2) $QN$ is $\mathrm{K}$-$\mathbb{P}_{t}$-subnormal in $G$.

If $q = p$, then $QN=Q$ is $\mathrm{K}$-$\mathbb{P}_{t}$-subnormal in $G$.

Let $q\not= p$. Write $H/N = QF(G)/N$. By Lemma 2.3(1) $QN/N$ is $\mathrm{K}$-$\mathbb{P}_{t}$-subnormal in $H/N$. Let's consider two cases.

1. Assume that $QN/N$ is subnormal in $H/N$. Since $QN/N$ is pronormal in $H/N$, by \cite[Chap. A, Lemma 6.3(d)]{DH} $QN/N$ is normal in $H/N$. Thus, $QN$ is normal in $H$. Since $\Phi(F(G))$ $\mathrm{char}$ $F(G)\unlhd G$, we have $\Phi(F(G))\unlhd G$. Then $N\leq \Phi(F(G))$. From $F(G)\unlhd H$ and by \cite[Chap. A, Theorem 9.2(e)]{DH} it follows that $\Phi(F(G))\leq \Phi(H)$. Thus, $N\leq \Phi(F(G))\leq \Phi(H)$. By Frattini argument $H=N_{H}(Q)QN=N_{H}(Q)$. Then we have $Q\unlhd QN$ and $Q$ is $\mathrm{K}$-$\mathbb{P}_{t}$-subnormal in $G$.

2. Assume that $QN/N$ is not subnormal in $H/N$. From $N\leq \Phi(G) < F(G)$ and $Q\in\mathrm{Syl}_{q}(H)$ we have $QN/N\not=H/N$. By Definition 1 it follow that there exists a chain of subgroups $QN/N\ = R_{0}/N \leq R_{1}/N \leq \cdots \leq R_{m-1}/N \leq R_{m}/N = H/N$
such that either $R_{i-1}/N\unlhd R_{i}/N$ or $|R_{i}/N : R_{i-1}/N|=p$ and $p-1$ is not divisible by the $(t+1)$th powers of primes for every $i = 1,\ldots , m$. Since $G/\Phi(G)\in \mathfrak{H}\subseteq \overline{\mathrm{w}}\mathfrak{U}$ and by Lemma 1.2 $\overline{\mathrm{w}}\mathfrak{U}$ is a hereditary saturated formation, we have $G\in\overline{\mathrm{w}}\mathfrak{U}$. Then $Q$ is $\mathrm{K}$-$\mathbb{P}$-subnormal in $G$ and $Q$ is $\mathrm{K}$-$\mathbb{P}$-subnormal in $QN$. Since $|QN:Q|=|N|$ and $p-1$ is not divisible by the $(t+1)$th powers of primes, this implies that $Q$ is $\mathrm{K}$-$\mathbb{P}_{t}$-subnormal in $G$.
Therefore, $G\in \mathfrak{H}$.
We arrive at a contradiction to the choice of $G$. Statement (6) has been proven.

Statement (7) follows from (2), (3) and (5).
%This completes the proof of the theorem.
\hspace{\stretch{1}}$\square$
%{\qedsymbol}
%\end{proof}

\medskip

In the work \cite{VasVasVeg16}, local definitions of the formation of groups whose Sylow subgroups are $\mathfrak{F}$-subnormal ($\mathrm{K}$-$\mathfrak{F}$-subnormal, respectively) were studied. Next we solve a similar problem for the case of $\mathrm{K}$-$\mathbb{P}_{t}$-subnormal Sylow subgroups.

\medskip

Since $\mathfrak{N}_{p}\mathfrak{A}(p-1)$ is a hereditary formation, the following result is easily verified.

\medskip

{\bf Lemma 3.2.} {\it Let $p$ be a prime number and $p-1$ is not divisible by the $(t+1)$th powers of prime numbers. The class of groups $(G\ | \ \mathrm{Syl}(G)\subseteq\mathfrak{N}_{p}\mathfrak{A}(p-1))$ is a hereditary formation.
}

\medskip

{\bf Теорема 3.2.}
{\it Let $\mathfrak{H}$ be the class of all groups in which every Sylow subgroup is $\mathrm{K}$-$\mathbb{P}_{t}$-subnormal. Then $\mathfrak{H}$  is a hereditary saturated formation that is defined by a local function $F$ such that
$F(p)=(G\in\mathfrak{S}\ | \ \mathrm{Syl}(G)\subseteq\mathfrak{N}_{p}\mathfrak{A}(p-1))$ if $p-1$ is not divisible by the $(t+1)$th powers of prime numbers;
$F(p)=\mathfrak{N}_{p}$ if $p-1$ is divisible by the $(t+1)$th power of some prime number.
}

\medskip

{\it Proof.}
%\begin{proof}
By Theorem 3.1 $\mathfrak{H}$ is a hereditary saturated formation.
Since $F(p)=(G\in\mathfrak{S}\ | \ \mathrm{Syl}(G)\subseteq\mathfrak{N}_{p}\mathfrak{A}(p-1))$ is a formation and $F(p)=\mathfrak{N}_{p}$ is a formation, $F$ is a local function. Let $\mathfrak{F}= LF(F)$. By \cite[Chap. IV, Proposition 3.14 and Theorem 4.6]{DH} $\mathfrak{F}$ is a hereditary saturated formation.

Show that $\mathfrak{F}\subseteq \mathfrak{H}$. Let $G$ be a group of least order in $\mathfrak{F}\setminus\mathfrak{H}$. Let $N$ be a minimal normal subgroup of $G$. Since $\mathfrak{F}\subseteq\mathfrak{S}$, $N$ is an $p$-group for some prime $p$. From $\mathfrak{N}_{p}\subseteq \mathfrak{F}\cap\mathfrak{H}$ it follows that $G\not= N$. We have $G/N\in\mathfrak{H}$, $\Phi(G)=1$ and $N$ is a unique minimal normal subgroup of $G$, since $\mathfrak{F}$ and $\mathfrak{H}$ are saturated formation. Then $G = NM$, where $M$ is a maximal subgroup of $G$, $\mathrm{Core}_{G}(M)=1$, $M\cap N=1$ and $N = C_{G}(N)$. Let $Q\in\mathrm{Syl}_{q}(G)$. Then $QN/N\in\mathrm{Syl}_{q}(G/N)$. From $G/N\in\mathfrak{H}$  by Lemma 2.2(3) we have that $QN$ is $\mathrm{K}$-$\mathbb{P}_{t}$-subnormal in $G$.

If $q=p$ then $QN=Q$ is $\mathrm{K}$-$\mathbb{P}_{t}$-subnormal in $G$.

Let $q\not=p$.
If $G\not=QN$, then $QN\in \mathfrak{H}$ since $\mathfrak{F}$ and $\mathfrak{H}$ are hereditary. Thus $Q$ is $\mathrm{K}$-$\mathbb{P}_{t}$-subnormal in $QN$ and $Q$ is $\mathrm{K}$-$\mathbb{P}_{t}$-subnormal in $G$.
Now suppose that $G=QN$. Then $Q=M$. We have that $Q\cong G/N=G/C_{G}(N)\in F(p)$ because  $G/N\in\mathfrak{F}$. Hence $F(p)\not=\mathfrak{N}_{p}$. Therefore $p-1$ is not divisible by the $(t+1)$th powers of prime numbers and $Q\in \mathfrak{A}(p-1)$. Then $G$ is supersoluble and $|N|=p$. We conclude that $|G:Q|=p$ and by Definition 1, we have that $Q$ is $\mathrm{K}$-$\mathbb{P}_{t}$-subnormal in $G$. Hence $G\in\mathfrak{H}$. It is a contradiction to the choice of $G$. Thus, $\mathfrak{F}\subseteq \mathfrak{H}$.

Prove that $\mathfrak{H}\subseteq\mathfrak{F}$.
Let $G$ be a group of least order in $\mathfrak{H}\setminus\mathfrak{F}$. From $G\in\mathfrak{H}$, we have that $G$ is soluble. We denote by $N$ a minimal normal subgroup of $G$. If $G=N$ then $G/C_{G}(N)=G/N\cong 1\in \mathfrak{F}$. It is a contradiction to the choice of $G$. Hence $G\not=N$.
Since $\mathfrak{H}$ and $\mathfrak{F}$ are saturated formations, $\Phi(G) = 1$. In $G$, $N$ is the unique minimal normal subgroup, $N = C_{G}(N) = F(G)$ and $|N| = p^{\alpha}$ for some prime $p$. The group $G= NM$, where $M$ is a maximal subgroup of $G$, $N\cap M = 1$ and $\mathrm{Core}_{G}(M)=1$. Since $G$ is Ore dispersive, we have $N\leq P\unlhd G$, $P\in\mathrm{Syl}_{p}(G)$ and $p$ is the largest prime number in $\pi(G)$. By \cite[cp. A, Theorem 15.6(b)]{DH} it follows that $O_{p}(M)=1$. From $P\cap M\leq O_{p}(M)$ we have that $P=(M\cap P)N=N\in\mathrm{Syl}_{p}(G)$ and $M$ is a $p'$-group.

Let $R\in\mathrm{Syl}_{q}(M)$. Then $R\in\mathrm{Syl}_{q}(G)$.

Suppose that $G\not= RN$. Let $H=RN$. Since $\mathfrak{H}$ is hereditary, we have that $H\in \mathfrak{H}$ and $H\in\mathfrak{F}$. Note that $C_{H}(N)=N$. Then $R\cong H/C_{H}(N)\in F(p)$. From $q\not=p$ it follows that $F(p)\not=\mathfrak{N}_{p}$ and  $p-1$ is not divisible by the $(t+1)$th powers of prime numbers. Thus $R\in \mathfrak{A}(p-1)$. Since $R$ is chosen arbitrarily, we have that $G/C_{G}(N)=G/N\cong M\in F(p)$. Hence $G\in\mathfrak{F}$. We have a contradiction to the choice of $G$.

Now suppose that $G= RN$. By \cite{Vas2004} $G$ is supersoluble. Then $|N|=p$. Hence $R=M\cong G/C_{G}(N)$ is isomorphic to a subgroup from $Z_{p-1}$. On the other hand, $R$ is $\mathrm{K}$-$\mathbb{P}_{t}$-subnormal in $G$, since $G\in\mathfrak{H}$. From $|G:R|=p$ it follows that $p-1$ is not divisible by the $(t+1)$th powers of prime numbers. Therefore $G/C_{G}(N)\in F(p)$ and $G\in\mathfrak{F}$. This contradicts the choice of $G$. Thus $\mathfrak{H}\subseteq\mathfrak{F}$. This means that the equality $\mathfrak{H}=\mathfrak{F}$ is proven.
\hspace{\stretch{1}}$\square$

%{\qedsymbol}
%\end{proof}

\medskip

{\bf Definition 3.1.} {\it  Denote by $\mathfrak{U}_{t}^{0}$ the class of all supersoluble groups in which every Sylow subgroup is $\mathrm{K}$-$\mathbb{P}_{t}$-subnormal.}

\medskip

{\bf Theorem 3.3.}
{\it The class of groups $\mathfrak{U}_{t}^{0}$ is a hereditary saturated formation that is defined by a local function $X$ such that
$X(p)=\mathfrak{N}_{p}\mathfrak{A}(p-1)$ if $p-1$ is not divisible by the $(t+1)$th powers of prime numbers;
$X(p)=\mathfrak{N}_{p}$ if $p-1$ is divisible by the $(t+1)$th power of some prime number.}

\medskip

{\it Proof.}
%\begin{proof}
We have that $\mathfrak{U}_{t}^{0}=\mathfrak{U}\cap\mathfrak{H}$.

Let $G\in\mathfrak{U}_{t}^{0}$ and let $H/K$ be its any chief factor. Then $|H/K|=p$ for some prime $p$. Suppose that $p-1$ is divisible by the $(t+1)$th power of some prime number. From $G\in\mathfrak{H}$ by Theorem 3.2, we have that $G/C_{G}(H/K)\in F(p)=\mathfrak{N}_{p}=X(p)$. Now suppose that $p-1$ is divisible by the $(t+1)$th power of some prime number. Since $G\in\mathfrak{U}$ we have that $G/C_{G}(H/K)\in\mathfrak{A}(p-1)\subseteq\mathfrak{N}_{p}\mathfrak{A}(p-1)=X(p)$. This means that $G\in LF(X)$, i.e. $\mathfrak{U}_{t}^{0}\subseteq LF(X)$.

Now let $G\in LF(X)$. Since $X(p)\subseteq\mathfrak{S}$ we have that $LF(X)\subseteq\mathfrak{S}$. Let $H/K$ be any chief factor of $G$. Then $H/K$ is abelian and $|H/K|=p^{\alpha}$ for some prime $p$.

Suppose that $p-1$ is divisible by the $(t+1)$th power of some prime number. Then $G/C_{G}(H/K)\in X(p)=\mathfrak{N}_{p}$.  By Theorem 3.2 $G\in\mathfrak{H}$. On the other hand, by \cite[Chap. A, Lemma 13.6(b)]{DH} $O_{p}(G/C_{G}(H/K))=1$ it follows that $G/C_{G}(H/K)=1\in\mathfrak{A}(p-1)$. Therefore $G\in\mathfrak{U}$. Thus $G\in\mathfrak{U}\cap\mathfrak{H}$.

Let $p-1$ is not divisible by the $(t+1)$th powers of prime numbers. Then $G/C_{G}(H/K)\in X(p)=\mathfrak{N}_{p}\mathfrak{A}(p-1)$. We have that $G/C_{G}(H/K)\in F(p)$ and $G\in\mathfrak{H}$ by Theorem 3.2. From \cite[Chap. A, Lemma 13.6(b)]{DH} we conclude that $G/C_{G}(H/K)\in \mathfrak{A}(p-1)$. Thus $G\in\mathfrak{U}$.  Therefore $G\in\mathfrak{U}_{t}^{0}$ and $LF(X)\subseteq \mathfrak{U}_{t}^{0}$.

Thus $LF(X)=\mathfrak{U}_{t}^{0}$.
\hspace{\stretch{1}}$\square$
%{\qedsymbol}
%\end{proof}

\medskip

{\bf Remark 3.1.} If $t=1$ then from Example 2.1 it follows that $\mathfrak{U}_{1}^{0}\subseteq s\mathfrak{U}$ and $\mathfrak{U}_{1}^{0}\not= s\mathfrak{U}$. Here $s\mathfrak{U}$ is a class of all  supersoluble groups in which each Sylow subgroup is submodular. The properties of submodular subgroups and the class $s\mathfrak{U}$ were studied in \cite{Zim} and \cite{VVA2015}.

\medskip

{\bf Теорема 3.4.}
{\it The class $\mathfrak{H}$ of all groups in which every Sylow subgroup is $\mathrm{K}$-$\mathbb{P}_{t}$-subnormal coincides with the class of all groups in which every Sylow subgroup is $\mathfrak{U}_{t}^{0}$-subnormal, i.e. $\mathfrak{H}=\mathrm{w}\mathfrak{U}_{t}^{0}$.
}

\medskip

{\it Proof.}
%\begin{proof}
From $\mathfrak{N}\subseteq \mathfrak{U}_{t}^{0}$ we have that $\pi(\mathfrak{U}_{t}^{0})=\mathbb{P}$.
Since $\mathfrak{U}_{t}^{0}$ is a hereditary saturated formation, by Lemma 1.5 $\mathrm{w}\mathfrak{U}_{t}^{0}$ is a hereditary saturated formation and $\mathrm{w}\mathfrak{U}_{t}^{0}\subseteq \mathfrak{S}$.

Show that $\mathfrak{H}\subseteq\mathrm{w}\mathfrak{U}_{t}^{0}$. Suppose that $\mathfrak{H}\setminus\mathrm{w}\mathfrak{U}_{t}^{0}\not=\varnothing$. Let $G$ be a group of the least order in $\mathfrak{H}\setminus\mathrm{w}\mathfrak{U}_{t}^{0}$. Then $G$ has the unique minimal normal subgroup $N$ and $\Phi(G) = 1$, since $\mathfrak{H}$ and $\mathrm{w}\mathfrak{U}_{t}^{0}$ are saturated formations. From $G\in\mathfrak{H}\subseteq\mathfrak{S}$ it follows that $N$ is a $p$-group for some prime $p$. Then $G=NM$, where $M$ is some maximal in $G$ subgroup with $\mathrm{Core}_{G}(M)=1$, $N\cap M=1$ and $N=C_{G}(N)$.

Let $Q$ be an arbitrary Sylow $q$-subgroup of $G$. Since $G/N\in\mathrm{w}\mathfrak{U}_{t}^{0}$ we have that $QN/N$ is $\mathfrak{U}_{t}^{0}$-subnormal in $G/N$. By Lemma 1.4(2) $QN$ is $\mathfrak{U}_{t}^{0}$-subnormal in $G$. If $q=p$, then $QN=Q$ is $\mathfrak{U}_{t}^{0}$-subnormal in $G$.

Let $q\not=p$. 
Suppose that $QN \not= G$. From $QN\in\mathrm{w}\mathfrak{U}_{t}^{0}$ we conclude that $Q$ is $\mathfrak{U}_{t}^{0}$-subnormal in $QN$. By Lemma 1.4(1) $Q$ is $\mathfrak{U}_{t}^{0}$-subnormal in $G$.
Let $QN = G$. Then $Q=M$. By Theorem 3.2 $Q\cong G/C_{G}(N) \in F(p)$. Since $G\not\in\mathfrak{N}$, we have that $F(p)=(G\in\mathfrak{S}\ | \ \mathrm{Syl}(G)\subseteq\mathfrak{N}_{p}\mathfrak{A}(p-1))$ if $p-1$ is not divisible by the $(t+1)$th powers of prime numbers. Hence $Q\in\mathfrak{A}(p-1)$ and $G\in\mathfrak{U}$. Thus $G\in \mathrm{w}\mathfrak{U}_{t}^{0}$ and $Q$ is $\mathfrak{U}_{t}^{0}$-subnormal in $G$. Consequently, $G\in\mathrm{w}\mathfrak{U}_{t}^{0}$. We get a contradiction. Hence $\mathfrak{H}\subseteq\mathrm{w}\mathfrak{U}_{t}^{0}$.

Prove that $\mathrm{w}\mathfrak{U}_{t}^{0}\subseteq\mathfrak{H}$. Suppose that $\mathrm{w}\mathfrak{U}_{t}^{0}\setminus\mathfrak{H}\not=\varnothing$. Choose a group $G$ of the least order in $\mathrm{w}\mathfrak{U}_{t}^{0}\setminus\mathfrak{H}$. It is clear that $G$ has the unique minimal normal subgroup $N = G^{\mathfrak{H}} = F(G)$ and $\Phi(G) = 1$. From $\mathfrak{H}\subseteq\mathfrak{S}$ it follows that $|N|=p^{\alpha}$ for some prime $p$. Let $S\in\mathrm{Syl}_{q}(G)$. By choice of $G$ we have $G\not=N$. From $G/N\in\mathfrak{H}$ and $SN/N\in\mathrm{Syl}_{q}(G/N)$ it follows that $SN/N$ is $\mathrm{K}$-$\mathbb{P}_{t}$-subnormal in $G/N$. By Lemma 2.2(2) $SN$ is $\mathrm{K}$-$\mathbb{P}_{t}$-subnormal in $G$. 

If $q=p$, then $SN=S$ is $\mathrm{K}$-$\mathbb{P}_{t}$-subnormal in $G$.

Suppose that $q\not=p$. 
If $SN \not= G$ then $SN\in\mathfrak{H}$ by choice of $G$. Consequently $S$ is $\mathrm{K}$-$\mathbb{P}_{t}$-subnormal in $SN$. By Lemma 2.1(2) $S$ is $\mathrm{K}$-$\mathbb{P}_{t}$-subnormal in $G$. 
Let $SN= G$. Then $S$ is maximal in $G$ and $\mathrm{Core}_{G}(S)=1$. From $G\in\mathrm{w}\mathfrak{U}_{t}^{0}$ it follows that $S$ is $\mathfrak{U}_{t}^{0}$-subnormal in $G$. Therefore $G^{\mathfrak{U}_{t}^{0}}\leq S$. Then $G^{\mathfrak{U}_{t}^{0}}\leq\mathrm{Core}_{G}(S)=1$ and $G\in\mathfrak{U}_{t}^{0}$. Thus $S$ is $\mathrm{K}$-$\mathbb{P}_{t}$-subnormal in $G$. Due to the arbitrary choice of $S$, we received a contradiction $G\in\mathfrak{H}$. Consequently $\mathrm{w}\mathfrak{U}_{t}^{0}\subseteq\mathfrak{H}$. Thus $\mathrm{w}\mathfrak{U}_{t}^{0}=\mathfrak{H}$.
\hspace{\stretch{1}}$\square$
%{\qedsymbol}
%\end{proof}

\medskip

This work was supported by the Ministry of Education of the Republic of Belarus (Grant no. 20211750"Convergence-2025") and by the Belarusian
Republican Foundation for Fundamental Research (Project F23RSF-237).

\renewcommand{\refname}{\hfill\normalsize\bf References\hfill}

%\noindent\begin{tabular}{m{8mm}m{100mm}}
%&{
A.\,F.\,Vasil'ev

Francisk Skorina Gomel State University,
Sovetskaya str., 104,
Gomel 246028, Belarus.
E-mail address: formation56amail.ru

\medskip

Т.\,I.\,Vasil'eva

Belarusian State University of Transport,
Kirov str., 34,
Gomel 246653, Belarus.
E-mail address: tivasilyeva@mail.ru

%\end{tabular}


\begin{thebibliography}{16}

\bibitem{VasVasTyu2010}
{\it Vasil'ev A.F., Vasil'eva T.I., Tyutyanov V.N.}
On the finite groups of supersoluble type // Siberian Math. J. 2010.
Vol.~51, No.~6. P. 1004--1012.

\bibitem{VasVasTyu2014}
{\it Vasil'ev A.F., Vasil'eva T.I., Tyutyanov V.N.}
On $\mathrm{K}-\mathbb{P}$-Subnormal Subgroups of Finite Groups //
Math. Notes. 2014. Vol. 95, No. 4. P. 471--480.

\bibitem{KniMon13}
{\it Kniahina V.N., Monakhov V.S.}
On supersolvability of finite groups with $\mathbb{P}$-subnormal subgroups // Internal. J. of Group Theory. 2013. Vol.2, No. 4. P. 21--29.

\bibitem{VasVasMys16}
{\it  Vasil'ev A.F., Vasil'eva T.I., Myslovets E.N.}
Finite widely $c$-supersoluble groups and their mutually permutable products // Siberian Math. J. 2016. Vol. 57, No. 3. pp. 476--485.

\bibitem{VasVasPar18}
{\it  Vasil'ev A.F., Vasil'eva T.I., Parfenkov K.L.}
Finite groups with three given subgroups //
Siberian Math. J. 2018. Vol. 59, No.1. P. 50--58.

\bibitem{BalLiPedSU20}
{\it Ballester-Bolinches A., Li Y., Pedraza-Aguilera M.C., Su N.}
On Factorised Finite Groups //
Mediterr. J. Math. 2020. Vol. 17, No. 2. P. 65.

%\bibitem{CheLiZha20} {\it Chen R. Li R. Zhao X.} Primary Subgroups and the Structure of Finite Groups // Commun. Algebr. 2020. Vol. 48, No. 10. P. 4436--4443.

\bibitem{Mur20}
{\it Murashka V.I.}
Finite Groups With Given Sets of $\mathfrak{F}$-Subnormal Subgroups // Asian-Eur. J. Math. 2020. Vol. 13, No. 4. P. 2050073.

\bibitem{LucNem21}
{\it Lucchini A., Nemmi D.}
The Non-$\mathfrak{F}$ Graph of a Finite Group //
Math. Nachr. 2021. Vol. 294, No. 10. P. 1912--1921.

\bibitem{Tro21}
{\it Trofimuk A.A.}
Finite factorizable groups with restrictions on factors.
Minsk: BSU Publishing Center, 2021 (In Russian).
%262 p.

\bibitem{VTI22}
{\it Vasilyeva T.I.}
Subgroups of the Fan of Sylow Subgroups and the Supersolvability of a Finite Group //
Math. Notes. 2021. Vol. 110, No. 2. P. 186--195.

\bibitem{BalMadShuPed22}
{\it Ballester-Bolinches A., Madanha S.Y., Shumba T.M.M., Pedraza-Aguilera M.C.}
On Certain Products of Permutable Subgroups //
Bull. Aust. Math. Soc. 2022. Vol. 105, No. 2. P. 278--285.

\bibitem{BalMadPedWU23}
{\it Ballester-Bolinches A., Madanha S.Y., Pedraza-Aguilera M.C., Wu X.}
On some products of finite groups //
%Proceedings of the Edinburgh Mathematical Society.
Proc. Edinburgh Math. Soc. 2023. Vol. 66, No. 1. P. 89.

\bibitem{CheZhaLi23}
{\it Chen R. Zhao X. Li X.}
$\mathbb{P}$-Subnormal Subgroups and the Structure of Finite Groups //
Ric. Mat. 2023. Vol. 72. P. 771--778.


\bibitem{VTIKor23}
{\it Vasilyeva T.I., Koranchuk  A.G.}
On Finite Groups with $\mathbb{P}_{\pi}$-Subnormal Subgroups //
Math. Notes. 2023. Vol. 114, No. 4. P. 421--432.

\bibitem{Mur24}
{\it Murashka V.I.}
Formations of finite groups in polynomial time: $\mathfrak{F}$-residuals and $\mathfrak{F}$-subnormality //
Journal of Symbolic Computation. 2024. Vol. 122. P. 102271.

\bibitem{YiXuKam24}
{\it Yi X., Xu Z., Kamornikov S.F.}
Finite groups with $\mathbb{P}$-subnormal Schmidt subgroups //
Trudy Instituta Matematiki i Mekhaniki UrO RAN. 2024. Vol. 30, No. 1. 100--108.



\bibitem{Lis24}
{\it Lisi F.}
A Jordan–Holder type theorem for finite groups //
Annali di Matematica. 2024. https://doi.org/10.1007/s10231-024-01456-w




\bibitem{BalEzq}
{\it Ballester-Bolinches A., Ezquerro L.M.}
Classes of Finite Groups, in Math. Appl. (Springer)
Dordrecht: Springer, 2006. Vol. 584.



\bibitem{MonSok_Sib_el23}
{\it Monakhov V.S., Sokhor I.L.}
Finite groups with formational subnormal primapy subgroups of bounded exponent // Siberian Electronic Mathematical News. 2023.  Vol. 20, No. 2.  P. 785--796.

\bibitem{DH} {\it Doerk K., Hawkes T.} Finite Soluble Groups. Berlin-New York:
Walter de Gruyter, 1992.  %898~p.

\bibitem{Wei}
Between Nilpotent and Solvable / H.G. Bray [and others]; edited by M. Weinstein.
Passaic: Polugonal Publishing House, 1982. %240 p.



\bibitem{Vas2004}
{\it Vasil'ev A.F.}
New properties of finite dinilpotent groups //
Vestsi Nats. Akad. Navuk Belarusi. Ser. Fiz.-Mat.
Navuk. 2004. No. 2. P. 29--33 (In Russian).

%\bibitem{Hup} {\it Huppert B.} Endliche Gruppen. I. / B. Huppert. Berlin: Springer, 1967. %795 s.

%\bibitem{sch1} {\it Schmidt R.}Modulare Untergruppen endlicher Gruppen / R.~Schmidt // J. Ill. Math. 1969. Vol. 13. P.~358--377.

\bibitem{VasVas11}
{\it Vasil'ev A.F., Vasil'eva T.I.}
On finite groups with generally subnormal Sylov subgroups //
Problems of physics, mathematics and technics. 2011. No. 4 (9). P. 86--91 (In Russian).

\bibitem{VasVasVeg16}
{\it Vasil'ev A.F., Vasil'eva T.I., Vegera A.S.}
Finite groups with generalized subnormal embedding of Sylov subgroups // Siberian Math. J. 2016. Vol. 57, No. 2. P. 200--212.

\bibitem{Zim}
{\it Zimmermann I.} Submodular subgroups in finite groups // Math. Z. 1989. Vol. 202. P.~545--557.

\bibitem{VVA2015}
{\it Vasilyev V.A.}
Finite groups with submodular Sylow subgroups // Siberian Math. J. 2015.
Vol.~56, No.~6. P. 1019--1027.




%\bibitem{Shem} %Шеметков, Л.А. Формации конечных групп / Л.А. Шеметков. %М.: Наука, 1978. %278~с.

%{\it Shemetkov L.A.} Formations of finite groups. Moscow: Nauka, 1987. (In Russian)

%Васильев А.Ф. О конечных группах сверхразрешимого типа /
%А.Ф. Васильев, Т.И. Васильева, В.Н. Тютянов
%// Сиб. мат. журн. --- 2010. --- Т.~51, №~6. --- С.~1270--1281.

%\bibitem{VasVasTju1} %Васильев, А.Ф. {\it Vasil'ev A.F., Vasil'eva T.I., Tyutyanov V.N.} On the products of $\Bbb{P}$-subnormal subgroups of finite groups // Siberian Math. J. 2012. Vol.~53, No.~1. P. 47--54.











\end{thebibliography}
\end{document}